\documentclass{amsart}

\usepackage{amsmath,amssymb,amsthm,caption,graphicx,interval,mathrsfs,mathtools,microtype,multicol,orcidlink,pgfplots,tikz,tikz-cd,xcolor}
\usepackage{ulem}
\usepackage{blkarray}

\renewcommand{\Re}{\mathrm{Re}}

\newcommand{\overbar}[1]{\mkern1.5mu \overline{\mkern-1.5mu #1 \mkern-1.5mu} \mkern1.5mu}

\newtheorem{counter}{}[section]

\theoremstyle{definition}

\newtheorem{remark}[counter]{Remark}

\theoremstyle{plain}
\newtheorem{corollary}[counter]{Corollary}
\newtheorem{lemma}[counter]{Lemma}
\newtheorem{proposition}[counter]{Proposition}
\newtheorem{theorem}[counter]{Theorem}

\title{A sharp bound for the functional calculus of $\rho$-contractions}
\date{}

\begin{document}

\thanks{The second named author is financed by the Dutch Research Council (NWO) grant OCENW.M20.292.}

\author[F.L.~Schwenninger]{Felix~L.~Schwenninger\,\orcidlink{0000-0002-2030-6504}}
\address{
University of Twente\\
P.O.~Box 217\\
7500~AE Enschede\\
The Netherlands}
\email{f.l.schwenninger@utwente.nl}

\author[J.~de~Vries]{Jens~de~Vries\,\orcidlink{0009-0008-2122-814X}}
\address{
University of Twente\\
P.O.~Box 217\\
7500~AE Enschede\\
The Netherlands}
\email{j.devries-4@utwente.nl}

\maketitle 

\begin{abstract}
Let $A$ be a $\rho$-contraction and $f$ a rational function mapping the closed unit disk into itself. With a new characterization of $\rho$-contractions we prove that
\begin{align*}
\big\|f(A)\big\|\leq \frac{\rho}{2}\big(1-|f(0)|^{2}\big)+\sqrt{\frac{\rho^{2}}{4}\big(1-|f(0)|^{2}\big){}^{2}+|f(0)|^{2}}.
\end{align*}
We further show that this bound is sharp. This refines an estimate by Okubo--Ando and, for $\rho=2$, is consistent with a result by Drury.
\end{abstract}
\section{Introduction}
Let $A$ be a bounded operator on a Hilbert space $H$. Suppose that $\rho\geq1$. A \textit{unitary $\rho$-dilation} of $A$ is a unitary $U$ on a Hilbert space $\mathfrak{H}$ that contains $H$ as a closed subspace such that 
\begin{align*}
A^{n}x = \rho\cdot\textbf{P}U^{n}x
\end{align*}
for all $x\in H$ and $n=1,2,\ldots$, where $\textbf{P}$ is the orthogonal projection from $\mathfrak{H}$ onto $H$. In 1966 Sz.-Nagy--Foia\c{s} \cite{nagy1966certain} proved that $A$ has a unitary $\rho$-dilation if and only if
\begin{align}\label{SzNF}
    1-2\Big(1-\frac{1}{\rho}\Big)\Re(\sigma^{*}A)-\Big(\frac{2}{\rho}-1\Big)A^{*}A\geq0
\end{align}
for all $\sigma\in\partial D$, where $D$ is the open unit disk in $\mathbb{C}$ centered at the origin. In 1968 Holbrook \cite{holbrook1968power} and Williams \cite{williams1968schwarz} introduced the \textit{numerical $\rho$-radius}
\begin{align*}
w_{\rho}(A)\coloneqq\inf\Big\{a>0:\frac{1}{a}A \ \text{has a unitary $\rho$-dilation}\Big\}
\end{align*}
of $A$ and argued that $A$ has a unitary $\rho$-dilation if and only if
\begin{align*}
    w_{\rho}(A)\leq1.
\end{align*}
Due to this fact, any $A$ that admits a unitary $\rho$-dilation is called a \textit{$\rho$-contraction}. It is well-known that $w_{1}(A)=\|A\|$ is the operator norm and $w_{2}(A)=w(A)$ is the standard numerical radius. Moreover, the numerical $\rho$-radius is bounded below by the spectral radius. For other key properties of the numerical $\rho$-radius we refer to \cite{holbrook1971inequalities} and \cite{ando1974convexity}. A complete description of the $\rho$-contractions on $\mathbb{C}^{2}$ was given by Okubo--Spitkovsky \cite{okubo2001characterization}.

Offering a different perspective, Pagacz--Pietrzycki--Wojtylak \cite{pagacz2020between} introduced the deformed numerical range, a convex set located between the spectrum of $A$ and the closed disk of radius $\|A\|$ centered at $0$, whose radius is exactly $w_{\rho}(A)$. This geometric construction agrees with the standard numerical range if $\rho=2$.

It follows from a result published in 1975 by Okubo--Ando \cite{okubo1975constants} that, if $A$ is a $\rho$-contraction, then 
\begin{align*}
    \|f(A)\|\leq\rho
\end{align*}
for all rational functions $f\colon\overbar{D}\to\mathbb{C}$ with $\|f\|_{\infty}\leq1$. In 2008 Drury \cite{drury2008symbolic} was able to refine Okubo--Ando's bound in the case $\rho=2$ by showing that, if $w(A)\leq1$, then
\begin{align}\label{k2bound}
\|f(A)\|\leq k\big(|f(0)|\big)
\end{align}
for all rational functions $f\in\overbar{D}\to\mathbb{C}$ with $\|f\|_{\infty}\leq1$, where $k\colon[0,1]\to[1,2]$ is defined by
\begin{align*}
k(s)\coloneqq\sqrt{2-3s^{2}+2s^{4}+2(1-s^{2})\sqrt{1-s^{2}+s^{4}}}
\end{align*}
for $s\in[0,1]$. 

In this article we refine Okubo--Ando's result for arbitrary $\rho\geq1$, see Theorem \ref{rhoDrury}. More precisely, we prove that, if $w_{\rho}(A)\leq1$, then
\begin{align*}
\|f(A)\|\leq k_{\rho}\big(|f(0)|\big)
\end{align*}
for all rational $f\colon\overbar{D}\to\mathbb{C}$ with $\|f\|_{\infty}\leq1$, where $k_{\rho}\colon[0,1]\to[1,\rho]$ is given by
\begin{align*}
k_{\rho}(s)\coloneqq\frac{\rho}{2}(1-s^{2})+\sqrt{\frac{\rho^{2}}{4}(1-s^{2})^{2}+s^{2}}
\end{align*}
for $s\in[0,1]$. In particular, the cases $\rho=1$ and $\rho=2$ correspond to von Neumann's inequality and Drury's bound \eqref{k2bound}, respectively. Finally, we prove that the established bound is sharp, which had not been stated by Drury even for $\rho=2$.

Drury's proof of \eqref{k2bound} relies on a result by Berger--Stampfli \cite[Theorem 4]{berger1967mapping} to reduce to the case where the rational function is a Blaschke factor. In the same spirit, we use a generalization for $\rho$-contractions of Berger--Stampfli's result to achieve a similar reduction, see Section \ref{Preservation} for the statement and references. This enables us to apply a new characterization of $\rho$-contractions, Theorem \ref{otherCharacterization}, as well as to carefully exploit the form of $k_{\rho}$ to prove the main result, Theorem \ref{rhoDrury}.

\section{Preliminaries}
Besides the relevant background material, this section presents a characterization of $\rho$-contractions, Theorem \ref{otherCharacterization}, which, to the best of our knowledge, appears to be new.

\subsection{Analytic functional calculi}
Given a compact subset $X$ of $\mathbb{C}$, write $\mathcal{C}(X)$ for the C*-algebra of continuous functions from $X$ to $\mathbb{C}$. Consider the disk algebra
\begin{align*}
    \mathcal{A}(\overbar{D})\coloneqq\{f\in\mathcal{C}(\overbar{D}):f \ \text{is analytic on $D$}\}.
\end{align*}
The rational functions on $\overbar{D}$ are uniformly dense in $\mathcal{A}(\overbar{D})$ by Mergelyan's theorem and the space $\mathcal{A}(\overbar{D})$ embeds isometrically into $\mathcal{C}(\partial D)$ by the maximum modulus principle.

Let $A$ be a bounded operator on a Hilbert space $H$. Further assume that the spectrum of $A$ is contained in the open disk $D$. The \textit{analytic functional calculus} $\gamma_{A}\colon\mathcal{A}(\overbar{D})\to\mathcal{L}(H)$ is defined by
\begin{align*}
    \gamma_{A}(f)\coloneqq\int_{\partial D}f(\sigma)(\sigma-A)^{-1} \ \frac{\mathrm{d}\sigma}{2\pi i}
\end{align*}
for $f\in\mathcal{A}(\overbar{D})$ and extends the rational functional calculus. For $\sigma\in\partial D$ we put
\begin{align*}
    P(\sigma,A)\coloneqq\frac{1}{\pi}\Re\big(\sigma(\sigma-A)^{-1}\big).
\end{align*}
It holds that
\begin{align}\label{integralRepresentationHFC}
\gamma_{A}(f)+f(0)=\int_{\partial D}f(\sigma)P(\sigma,A) \ |\mathrm{d}\sigma|
\end{align}
for all $f\in\mathcal{A}(\overbar{D})$, see for example \cite[Section I.1]{cassier1996contractions}. The integral representation \eqref{integralRepresentationHFC}, formulated here for the disk, is a special case of a more general formula that encompasses more general domains, see \cite{delyon1999generalization}, \cite{crouzeix2007numerical} and \cite{schwenninger2024double} and the references therein.

\subsection{Two characterizations of $\rho$-contractions}
Suppose that $\rho\geq1$. The following result, together with the short proof, can be found between the lines in \cite[Section II.1]{cassier1996contractions}.
\begin{proposition}\label{potentialCharacterization}
Let $A$ be a bounded operator on a Hilbert space $H$ such that $D$ contains the spectrum of $A$. Then $A$ is a $\rho$-contraction if and only if
\begin{align}\label{potentialCharacterizationEquation}
P(\sigma,A)\geq\frac{2-\rho}{2\pi}
\end{align}
for all $\sigma\in\partial D$.
\end{proposition}
\begin{proof}
    One readily verifies that
    \begin{align*}
    &2\pi P(\sigma,A)-(2-\rho)\\
    &\qquad=\sqrt{\rho}(\sigma^{*}-A^{*})^{-1}\Big(1-2\Big(1-\frac{1}{\rho}\Big)\Re(\sigma^{*}A)-\Big(\frac{2}{\rho}-1\Big)A^{*}A\Big)\sqrt{\rho}(\sigma-A)^{-1}
    \end{align*}
    for all $\sigma\in\partial D$. Hence \eqref{SzNF} holds for all $\sigma\in\partial D$ if and only if \eqref{potentialCharacterizationEquation} holds for all $\sigma\in\partial D$.
\end{proof}
Next we present a new characterization of $\rho$-contractions, which plays an essential role in the proof of the main result. The case $\rho=2$ aligns with a result by Klaja--Mashreghi--Ransford \cite[Theorem 3.2]{klaja2016mapping}, which was obtained via different techniques.
\begin{theorem}\label{otherCharacterization}
Let $A$ be a bounded operator on a Hilbert space $H$ such that $D$ contains the spectrum of $A$. Then $A$ is a $\rho$-contraction if and only if
\begin{align}\label{otherCharacterizationEquation}
	1-2\Big(1-\frac{1}{\rho}\Big)\Re(z^{*}A)-\Big(\frac{1}{\rho^{2}}-\Big(1-\frac{1}{\rho}\Big){}^{2}|z|^{2}\Big)A^{*}A\geq0
\end{align}
for all $z\in\overbar{D}$.
\end{theorem}
\begin{proof}
For each $z\in\overbar{D}$ we define the resolvent operator
	\begin{align*}
		S_{\rho}(z)\coloneqq\Big(1-\Big(\frac{1}{\rho}+\Big(1-\frac{1}{\rho}\Big)z^{*}\Big)A\Big){}^{-1},
	\end{align*}
	where we use that the spectrum of $A$ is contained in $D$. A laborious but elementary computation shows that
	\begin{align*}
        &S_{\rho}(z)^{*}\Big(1-2\Big(1-\frac{1}{\rho}\Big)\Re(z^{*}A)-\Big(\frac{1}{\rho^{2}}-\Big(1-\frac{1}{\rho}\Big){}^{2}|z|^{2}\Big)A^{*}A\Big)S_{\rho}(z)
        \\
        &\qquad=\Re\Big(1+\frac{2}{\rho}AS_{\rho}(z)\Big)
	\end{align*}
	for all $z\in\overbar{D}$. For any $x\in H$ the real-valued continuous function
 \begin{align*}
     z\mapsto \Big\langle\Re\Big(1+\frac{2}{\rho}AS_{\rho}(z)\Big)x,x\Big\rangle,\qquad z\in\overbar{D}
 \end{align*}
 is harmonic on $D$ and therefore attains its global minimum on $\partial D$. It follows that \eqref{SzNF} holds for all $\sigma\in\partial D$ if and only if \eqref{otherCharacterizationEquation} holds for all $z\in\overbar{D}$.
 \end{proof}
 
\section{Functional calculus of $\rho$-contractions}
Let $A$ be a bounded operator on a Hilbert space $H$. Suppose that $\rho\geq1$. If $A$ is a $\rho$-contraction on a Hilbert space $H$, then the spectrum of $A$ is contained in $\overbar{D}$. In particular, for any rational function $f\colon \overbar{D}\to\mathbb{C}$ with poles off $\overbar{D}$ the operator $f(A)$ is well-defined.

\subsection{Preservation of $\rho$-contractivity}\label{Preservation}
The following result, which states that $f(A)$ is again a $\rho$-contraction if $\|f\|_{\infty}\leq1$ and $f(0)=0$, has been obtained through different methods: Sz.-Nagy--Foia\c{s} \cite[Proposition 4]{nagy1966certain} used the unitary $\rho$-dilation of $A$ directly, while Williams \cite[Theorem 2]{williams1968schwarz} provided an argument based on techniques from Berger--Stampfli \cite[Theorem 4]{berger1967mapping}. Under the assumption that the spectrum of $A$ is contained in $D$, the idea of the latter approach is to find a positive operator-valued Borel measure $E$ on $\partial D$ such that
\begin{align}\label{integralRepresentation}
P\big(\tau,\gamma_{A}(f)\big)-\frac{2-\rho}{2\pi}=\lim_{\varepsilon\nearrow1}\int_{\partial D}\Re\Big(\frac{\tau+\varepsilon f(\sigma)}{\tau-\varepsilon f(\sigma)}\Big) \ \frac{\mathrm{d}E\sigma}{2\pi}
\end{align}
for all $\tau\in\partial D$ and exploit Proposition \ref{potentialCharacterization}. Using a representation theorem by Herglotz--Riesz, Williams showed the existence of $E$ implicitly. A different approach by Cassier--Fack \cite[Proposition 6]{cassier1996contractions} was based on an explicit expression for $E$, namely,
\begin{align}\label{operatorMeasure}
    E(\Gamma)\coloneqq \int_{\Gamma}P(\sigma,A) \ |\mathrm{d}\sigma|-\frac{2-\rho}{2\pi}|\Gamma|
\end{align}
for Borel sets $\Gamma\subseteq\partial D$. In the proof below we derive \eqref{integralRepresentation} with the explicit expression \eqref{operatorMeasure} as well, while adhering more closely to Berger--Stampfli's techniques.

\begin{theorem}\label{rhoBS}
    Let $A$ be a $\rho$-contraction on a Hilbert space $H$. Then $f(A)$ is a $\rho$-contraction on $H$ whenever $f\colon\overbar{D}\to\mathbb{C}$ is a rational function with $\|f\|_{\infty}\leq1$ and $f(0)=0$.
\end{theorem}
\begin{proof}
        Let $f\colon\overbar{D}\to\mathbb{C}$ be a rational function with $\|f\|_{\infty}\leq1$ and $f(0)=0$. 
        
        First consider the case where the spectrum of $A$ is contained in $D$. It follows that the spectrum of $\gamma_{A}(f)$ is also contained in $D$ by the maximum modulus principle. So by Proposition \ref{potentialCharacterization} it suffices to prove that 
        \begin{align*}
            P\big(\tau,\gamma_{A}(f)\big)\geq\frac{2-\rho}{2\pi}
        \end{align*}
        for all $\tau\in\partial D$. Let $\tau\in\partial D$ and $0<\varepsilon<1$ be arbitrary. We compute
	   \begin{align*}
            P\big(\tau,\gamma_{A}(\varepsilon f)\big)&=\frac{1}{\pi}\Re\big(\tau(\tau-\gamma_{A}(\varepsilon f))^{-1}\big)\\
	       &=\frac{1}{\pi}\Re\Big(1+\sum_{n=1}^{\infty}\frac{\varepsilon^{n}}{\tau^{n}}\gamma_{A}(f^{n})\Big)\\
	       &=\frac{1}{2\pi}\Re\Big(\gamma_{A}\Big(1+2\sum_{n=1}^{\infty}\frac{\varepsilon^{n}}{\tau^{n}}f^{n}\Big)+1\Big)\\
           &=\frac{1}{2\pi}\Re\big(\gamma_{A}\big((\tau+\varepsilon f)(\tau-\varepsilon f)^{-1}\big)+1\big)
	   \end{align*}
        Since $f(0)=0$, we infer from \eqref{integralRepresentationHFC} and Cauchy's integral formula that
        \begin{align*}
		 &P\big(\tau,\gamma_{A}(\varepsilon f)\big)-\frac{2-\rho}{2\pi}\\
         &\qquad=\frac{1}{2\pi}\Re\Big(\int_{\partial D}\frac{\tau+\varepsilon f(\sigma)}{\tau-\varepsilon f(\sigma)}P(\sigma,A) \ |\mathrm{d}\sigma|\Big)-\frac{1}{2\pi}\Re\Big(\int_{\partial D}\frac{\tau+\varepsilon f(\sigma)}{\tau-\varepsilon f(\sigma)} \ |\mathrm{d}\sigma|\Big)\frac{2-\rho}{2\pi}\\
          &\qquad=\int_{\partial D}\Re\Big(\frac{\tau+\varepsilon f(\sigma)}{\tau-\varepsilon f(\sigma)}\Big) \ \frac{\mathrm{d}E\sigma}{2\pi},
        \end{align*}
        where $E$ is the positive operator-valued Borel measure defined by \eqref{operatorMeasure}. In particular, we recover the formula \eqref{integralRepresentation}, which yields the desired estimate as the M\"obius transformation
    \begin{align*}
        z\mapsto \frac{\tau+\varepsilon z}{\tau-\varepsilon z},\qquad z\in\overbar{D}
    \end{align*}
    maps into $\{z\in\mathbb{C}:\Re(z)\geq0\}$.

    The case where the spectrum of $A$ intersects $\partial D$ follows from approximating $f(A)$ by $f(rA)$ for sufficiently large $0<r<1$.
	\end{proof}

When $\rho=2$, Theorem \ref{rhoBS} matches the celebrated numerical range mapping result by Berger--Stampfli \cite[Theorem 4]{berger1967mapping}. By slightly modifying the proof of Theorem \ref{rhoBS}, one obtains the more general numerical range mapping result by Putinar--Sandberg \cite[Theorem 3]{putinar2005skew} for convex domains beyond the disk, see \cite[Theorem 4.1]{schwenninger2024double}.

\begin{remark}
    Theorem \ref{rhoBS} can be reformulated as follows. If $w_{\rho}(A)\leq1$, then
\begin{align*}
    w_{\rho}\big(f(A)\big)\leq1
\end{align*} 
for every rational function $f\colon\overbar{D}\to\mathbb{C}$ with $\|f\|_{\infty}\leq1$ and $f(0)=0$. In the situation $f(0)\neq0$, a weaker bound remains valid. The case $\rho=2$ was studied by Drury \cite{drury2008symbolic} in 2008. In that work it was shown that, if $w(A)\leq1$, then 
\begin{align*}
    w\big(f(A)\big)\leq\frac{5}{4}
\end{align*}
for every rational function $f\colon\overbar{D}\to\mathbb{C}$ with $\|f\|_{\infty}\leq1$. The general case was later covered by Badea--Crouzeix--Klaja \cite{badea2018spectral} in 2018, see also \cite{davidson2018complete}. They proved that, if $w_{\rho}(A)\leq1$, then
\begin{align*}
w_{\rho}\big(f(A)\big)\leq\frac{\rho^{2}+1+(\rho-1)\sqrt{5\rho^{2}+2\rho+1}}{2\rho^{2}}
\end{align*}
for every rational function $f\colon\overbar{D}\to\mathbb{C}$ with $\|f\|_{\infty}\leq1$.
\end{remark}

\begin{remark}
Cassier--Suciu \cite[Theorem 5.1]{cassier2006mapping} showed that, if $A$ is a $\rho$-contraction, then $f(A)$ is a $\rho_{f}$-contraction for any non-constant rational function $f\colon\overbar{D}\to\mathbb{C}$ with $\|f\|_{\infty}\leq1$, where
\begin{align*}
    \rho_{f}\coloneqq1+(\rho-1)\frac{1+|f(0)|}{1-|f(0)|}.
\end{align*}
Note that $\rho_{f}=\rho$ if and only if $f(0)=0$.
\end{remark}

\subsection{Main result}
Consider the function $k_{\rho}\colon[0,1]\to[1,\rho]$ defined by
\begin{align}\label{k-function}
k_{\rho}(s)\coloneqq\frac{\rho}{2}(1-s^{2})+\sqrt{\frac{\rho^{2}}{4}(1-s^{2})^{2}+s^{2}}
\end{align}
for $s\in[0,1]$.

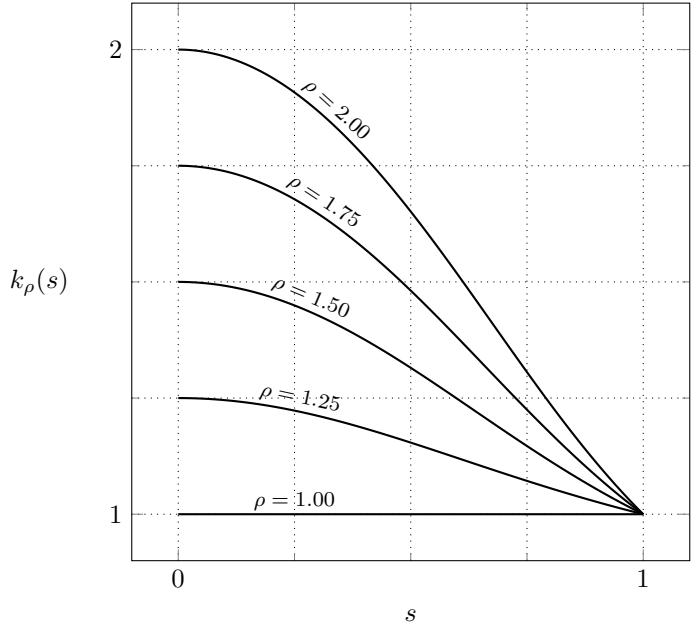
\begin{figure}[h]
	\begin{tikzpicture}[trim axis left]
		\begin{axis}[
			width=9cm,
			height=9cm,
			grid=both,
			grid style={dotted,black},
			xlabel={$s$},
			ylabel={$k_{\rho}(s)$},
			ylabel style={rotate=-90},
			xtick={0,1},
			minor x tick num=3,
			ytick={1,2},
			minor y tick num=3
			]
			\addplot[domain=0:1,smooth,samples=500,thick]
			{(1-x^2)+sqrt((1-x^2)^2+x^2)}
			node [sloped,yshift=5pt,font=\footnotesize,pos=0.25]	{$\rho=2.00$};
			\addplot[domain=0:1,smooth,samples=500,thick]
			{(1.75/2)*(1-x^2)+sqrt((1.75/2)^2*(1-x^2)^2+x^2)}
			node [sloped,yshift=5pt,font=\footnotesize,pos=0.25]	{$\rho=1.75$};
			\addplot[domain=0:1,smooth,samples=500,thick]
			{(1.5/2)*(1-x^2)+sqrt((1.5/2)^2*(1-x^2)^2+x^2)}
			node [sloped,yshift=5pt,font=\footnotesize,pos=0.25]	{$\rho=1.50$};
			\addplot[domain=0:1,smooth,samples=500,thick]
			{(1.25/2)*(1-x^2)+sqrt((1.25/2)^2*(1-x^2)^2+x^2)}
			node [sloped,yshift=5pt,font=\footnotesize,pos=0.25]	{$\rho=1.25$};
			\addplot[domain=0:1,smooth,samples=500,thick]
			{1}
			node [sloped,yshift=5pt,font=\footnotesize,pos=0.25]	{$\rho=1.00$};
		\end{axis}
	\end{tikzpicture}
    \caption{Graph of $k_{\rho}$ for several values of $\rho$.}
    \end{figure}

Note that $k_{\rho}(0)=\rho$ and $k_{\rho}(1)=1$. The function is constant for $\rho=1$ and strictly decreasing for $\rho\neq1$. Moreover, it is clear that \eqref{k-function} satisfies the quadratic relation
\begin{align}\label{quadraticRelation}
    k_{\rho}(s)^{2}-\rho(1-s^{2})k_{\rho}(s)-s^{2}=0,
\end{align}
for all $s\in[0,1]$. The following technical lemma will be used in the proof of the main result.

 \begin{lemma}\label{technicalLemma}
 For $\rho>1$ it holds that
 \begin{align*}
     -\Big(1-\frac{1}{\rho}\Big)<\frac{s-sk_{\rho}(s)^{2}}{k_{\rho}(s)^{2}-s^{2}}<0
 \end{align*}
 for all $0<s<1$.
 \end{lemma}
 \begin{proof}
It follows from \eqref{quadraticRelation} that
 \begin{align*}
     \frac{s-sk_{\rho}(s)^{2}}{k_{\rho}(s)^{2}-s^{2}}=-\Big(s-\frac{s}{\rho k_{\rho}(s)}\Big).
 \end{align*}
The derivative of the function $F_{\rho}\colon\interval[open]{0}{1}\to\mathbb{R}$ defined by
     \begin{align*}
         F_{\rho}(s)\coloneqq -\Big(s-\frac{s}{\rho k_{\rho}(s)}\Big)
     \end{align*}
    for $0<s<1$ satisfies
     \begin{align*}
     F_{\rho}'(s)=\frac{1}{2s^{2}}\cdot\frac{(1-s^{2})\big(k_{\rho}(s)-\rho\big)}{k_{\rho}(s)-\dfrac{\rho}{2}(1+s^{2})}<0.
     \end{align*}
     Thus $F_{\rho}$ is strictly decreasing and the claim follows.
 \end{proof}

We are now ready to prove the main result.

\begin{theorem}\label{rhoDrury}
    Let $A$ be a $\rho$-contraction on a Hilbert space $H$. Then
	\begin{align*}
		\|f(A)\|\leq k_{\rho}\big(|f(0)|\big)
	\end{align*}
    for every rational function $f\colon\overbar{D}\to\mathbb{C}$ with $\|f\|_{\infty}\leq1$.
\end{theorem}
\begin{proof}
Let $f\colon\overbar{D}\to\mathbb{C}$ be a rational function with $\|f\|_{\infty}\leq1$. For each $\lambda\in D$ we consider the function $g_{\lambda}\colon\overbar{D}\to\mathbb{C}$ given by
\begin{align*}
    g_{\lambda}(z)\coloneqq\frac{\lambda+z}{1+\lambda^{*}z}
\end{align*}
for $z\in\overbar{D}$.

First consider the case where the spectrum of $A$ is contained in $D$. By Theorem \ref{rhoBS} we may assume that $|f(0)|>0$. By the maximum modulus principle we may assume that $|f(0)|<1$. The function $h\colon\overbar{D}\to\mathbb{C}$ defined by
\begin{align*}
    h(z)\coloneqq g_{-f(0)}\big(f(z)\big)
\end{align*}
for $z\in\overbar{D}$ satisfies $\|h\|_{\infty}\leq1$ and $h(0)=0$. Fix $\zeta\in\partial D$ such that $\zeta f(0)=|f(0)|$. It holds that
\begin{align*}
    \zeta\gamma_{A}(f)=\zeta g_{f(0)}\big(\gamma_{A}(h)\big)=g_{\zeta f(0)}\big(\zeta\gamma_{A}(h)\big)
\end{align*}
and therefore $\|\gamma_{A}(f)\|=\|g_{s}(C)\|$ with $C\coloneqq\zeta\gamma_{A}(h)$ and $s\coloneqq|f(0)|$. Theorem $\ref{rhoBS}$ implies that $C$ is a $\rho$-contraction. The spectrum of $C$ is contained in $D$ by the maximum modulus principle and Lemma \ref{technicalLemma} asserts that
\begin{align*}
-\Big(1-\frac{1}{\rho}\Big)\leq\frac{s-sk_{\rho}(s)^{2}}{k_{\rho}(s)^{2}-s^{2}}\leq0.
\end{align*}
Hence Theorem \ref{otherCharacterization} implies that
\begin{align*}
1-\Big(\frac{s-sk_{\rho}(s)^{2}}{k_{\rho}(s)^{2}-s^{2}}\Big)(C^{*}+C)-\Big(\frac{1}{\rho^{2}}-\Big(\frac{s-sk_{\rho}(s)^{2}}{k_{\rho}(s)^{2}-s^{2}}\Big){}^{2}\Big)C^{*}C\geq0.
\end{align*}
Using \eqref{quadraticRelation} we quickly find that
\begin{align*}
\frac{1}{\rho^{2}}-\Big(\frac{s-sk_{\rho}(s)^{2}}{k_{\rho}(s)^{2}-s^{2}}\Big){}^{2}=\frac{1-s^{2}k_{\rho}(s)^{2}}{k_{\rho}(s)^{2}-s^{2}}.
\end{align*}
Consequently, we have
\begin{align*}
1-\Big(\frac{s-sk_{\rho}(s)^{2}}{k_{\rho}(s)^{2}-s^{2}}\Big)(C^{*}+C)-\Big(\frac{1-s^{2}k_{\rho}(s)^{2}}{k_{\rho}(s)^{2}-s^{2}}\Big)C^{*}C\geq0.
\end{align*}
After rescaling and rearraging this amounts to
\begin{align*}
    k_{\rho}(s)^{2}(1+sC)^{*}(1+sC)-(s+C)^{*}(s+C)\geq0,
\end{align*}
which is in turn equivalent to $\|g_{s}(C)\|\leq k_{\rho}(s)$.

The case where the spectrum of $A$ intersects $\partial D$ follows from approximating $f(A)$ by $f(rA)$ for sufficiently large $0<r<1$.
\end{proof}
We immediately recover Okubo--Ando's bound \cite[Theorem 2]{okubo1975constants}.
\begin{corollary}
Let $A$ be a $\rho$-contraction on a Hilbert space $H$. Then
\begin{align*}
    \|f(A)\|\leq \rho
\end{align*}
for every rational function $f\colon\overbar{D}\to\mathbb{C}$ with $\|f\|_{\infty}\leq1$.
\end{corollary}
\subsection{Sharpness}
The bound in Theorem \ref{rhoDrury} is sharp. To see this, consider for a fixed $\rho\geq1$ the matrix
\begin{align*}
        A\coloneqq\begin{bmatrix}
            0&\rho\\
            0&0
        \end{bmatrix}.
\end{align*}
It is well-known that $A$ is a $\rho$-contraction on $\mathbb{C}^{2}$, see for example \cite[Theorem 4.3]{holbrook1968power}. We claim that for every $0\leq s\leq1$ there exists a rational function $g_{s}\colon\overbar{D}\to\mathbb{C}$ with $\|g_{s}\|_{\infty}\leq1$ and $|g_{s}(0)|=s$ such that
\begin{align*}
\|g_{s}(A)\|=k_{\rho}(s).
\end{align*}
For any rational function $f\colon\overbar{D}\to\mathbb{C}$ it holds that
\begin{align*}
    f(A)=\begin{bmatrix}
        f(0)&\rho f'(0)\\
        0&f(0)
    \end{bmatrix}
\end{align*}
and therefore
\begin{align*}
    \|f(A)\|=\frac{\rho}{2}|f'(0)|+\sqrt{\frac{\rho^{2}}{4}|f'(0)|^{2}+|f(0)|^{2}}.
\end{align*}
For $s=1$ we take the constant function $g_{s}\coloneqq1$. For $s\neq1$ we define $g_{s}\colon\overbar{D}\to\mathbb{C}$ by
\begin{align*}
    g_{s}(z)\coloneqq\dfrac{s+z}{1+sz}
\end{align*}
for $z\in\overbar{D}$. Clearly, $\|g_{s}\|_{\infty}=1$ and $g_{s}(0)=s$. Moreover, $g_{s}'(0)=1-s^{2}$. This proves the desired equality.

\begin{remark}
The first part of the proof of Theorem \ref{rhoDrury} reduces the problem to the case where the rational function is a Blaschke factor. Hence it is not surprising that sharpness is attained at a Blaschke factor.
\end{remark}

\bibliographystyle{amsplain}
\bibliography{Bibliography}

\providecommand{\bysame}{\leavevmode\hbox to3em{\hrulefill}\thinspace}
\providecommand{\MR}{\relax\ifhmode\unskip\space\fi MR }
% \MRhref is called by the amsart/book/proc definition of \MR.
\providecommand{\MRhref}[2]{%
  \href{http://www.ams.org/mathscinet-getitem?mr=#1}{#2}
}
\providecommand{\href}[2]{#2}
\begin{thebibliography}{10}

\bibitem{ando1974convexity}
T.~Ando and K.~Nishio, \emph{Convexity properties of operator radii associated with unitary $\rho$-dilations.}, Michigan Mathematical Journal \textbf{20} (1974), no.~4, 303--307.

\bibitem{badea2018spectral}
C.~Badea, M.~Crouzeix, and H.~Klaja, \emph{Spectral sets and operator radii}, Bulletin of the London Mathematical Society \textbf{50} (2018), no.~6, 986--996.

\bibitem{berger1967mapping}
C.A. Berger and J.G. Stampfli, \emph{Mapping theorems for the numerical range}, American Journal of Mathematics \textbf{89} (1967), no.~4, 1047--1055.

\bibitem{cassier1996contractions}
G.~Cassier and T.~Fack, \emph{Contractions in von {N}eumann algebras}, Journal of Functional Analysis \textbf{135} (1996), no.~2, 297--338.

\bibitem{cassier2006mapping}
G.~Cassier and N.~Suciu, \emph{Mapping theorems and {H}arnack ordering for $\rho$-contractions}, Indiana University Mathematics Journal (2006), 483--523.

\bibitem{crouzeix2007numerical}
M.~Crouzeix, \emph{Numerical range and functional calculus in {H}ilbert space}, Journal of Functional Analysis \textbf{244} (2007), no.~2, 668--690.

\bibitem{davidson2018complete}
K.~Davidson, V.~Paulsen, and H.~Woerdeman, \emph{Complete spectral sets and numerical range}, Proceedings of the American Mathematical Society \textbf{146} (2018), no.~3, 1189--1195.

\bibitem{delyon1999generalization}
B.~Delyon and F.~Delyon, \emph{Generalization of von {N}eumann's spectral sets and integral representation of operators}, Bulletin de la Soci{\'e}t{\'e} Math{\'e}matique de France \textbf{127} (1999), no.~1, 25--41.

\bibitem{drury2008symbolic}
S.W. Drury, \emph{Symbolic calculus of operators with unit numerical radius}, Linear Algebra and its Applications \textbf{428} (2008), no.~8-9, 2061--2069.

\bibitem{holbrook1968power}
J.A.R. Holbrook, \emph{On the power-bounded operator of {S}z.-{N}agy and {F}oia\c{s}}, Acta Scientiarum Mathematicarum \textbf{29} (1968), 299--310.

\bibitem{holbrook1971inequalities}
\bysame, \emph{Inequalities governing the operator radii associated with unitary $\rho$-dilations.}, Michigan Mathematical Journal \textbf{18} (1971), no.~2, 149--159.

\bibitem{klaja2016mapping}
H.~Klaja, J.~Mashreghi, and T.~Ransford, \emph{On mapping theorems for numerical range}, Proceedings of the American Mathematical Society \textbf{144} (2016), no.~7, 3009--3018.

\bibitem{okubo1975constants}
K.~Okubo and T.~Ando, \emph{Constants related to operators of class ${C}_{\rho}$}, Manuscripta Mathematica \textbf{16} (1975), no.~4, 385--394.

\bibitem{okubo2001characterization}
K.~Okubo and I.~Spitkovsky, \emph{On the characterization of {$2\times 2$} {$\rho$}-contraction matrices}, Linear Algebra and its Applications \textbf{325} (2001), no.~1-3, 177--189.

\bibitem{pagacz2020between}
P.~Pagacz, P.~Pietrzycki, and M.~Wojtylak, \emph{Between the von {N}eumann inequality and the {C}rouzeix conjecture}, Linear Algebra and its Applications \textbf{605} (2020), 130--157.

\bibitem{putinar2005skew}
M.~Putinar and S.~Sandberg, \emph{A skew normal dilation on the numerical range of an operator}, Mathematische Annalen \textbf{331} (2005), no.~2, 345--357.

\bibitem{schwenninger2024double}
F.L. Schwenninger and J.~de~Vries, \emph{The double-layer potential for spectral constants revisited}, arXiv:2409.15954, preprint (2024).

\bibitem{nagy1966certain}
B.~Sz.-Nagy and C.~Foia\c{s}, \emph{On certain classes of power-bounded operators in {H}ilbert space}, Acta Scientiarum Mathematicarum \textbf{27} (1966), 17--25.

\bibitem{williams1968schwarz}
J.P. Williams, \emph{Schwarz norms for operators}, Pacific Journal of Mathematics \textbf{24} (1968), no.~1, 181--188.

\end{thebibliography}

\end{document}